\input amstex
\documentstyle{amsppt}

\magnification=1200
\parskip 5pt
\pagewidth{5.4in} \pageheight{7.1in}

\baselineskip=14pt \expandafter\redefine\csname
logo\string@\endcsname{} \NoBlackBoxes \NoRunningHeads

\redefine\qed{$\ \blacksquare$}

\def\r1{\sqrt{-1}}
\def\g{\frak{g}}

\def\s{\frak{s}}
\def\h{\frak{h}}
\def\r1{\sqrt{-1}}

\def\H{\hat{H}}
\def\quott{/\! /}
\def\bZ{\mathord{\Bbb Z}}
\def\bR{\mathord{\Bbb R}}

\def\bC{\mathord{\Bbb C}}

\topmatter
\title The kernel of the equivariant Kirwan map and the residue formula
\endtitle
\author{Lisa C. Jeffrey and Augustin-Liviu Mare  }
\endauthor
\abstract {Using the notion of equivariant Kirwan map, as defined
by Goldin \cite{Go}, we prove that --- in the case of Hamiltonian
torus actions with isolated fixed points --- Tolman and Weitsman's
description of the kernel of the Kirwan map  can be deduced
directly from the residue theorem of \cite{Je-Ki1} and
\cite{Je-Ki2}. A characterization of the kernel of the Kirwan
map in terms of residues of one variable (i.e.
 associated to Hamiltonian {\it circle} actions) is obtained.  }
\endabstract
\endtopmatter
\document

\head \S 1 Introduction
\endhead

 Let $H$ be a torus acting in a
Hamiltonian fashion on a compact symplectic manifold $N$ and $S\subset H$
a circular subgroup. Assume that $0\in \s^*$ is a regular value of
the moment map $\mu_S:M\to \s^*$ which corresponds to the $S$
action on $N$. Define the symplectic reduction
$N_{\text{red}}=\mu_S^{-1}(0)/S$ and consider
$$\kappa_S :H_H^*(N)\to H_H^*(\mu_S^{-1}(0))=H_{H/S}^*(N_{\text{red}}).$$
 R. Goldin \cite{Go}  called this map the {\it equivariant Kirwan map} and then
 she proved:

 \proclaim{Theorem 1.1}{\rm (see [Go, Theorem 1.2])} The map $\kappa_S$ is
 surjective.
 \endproclaim

Assume that the fixed point set ${\Cal F}=N^S$ consists of
isolated points. Consider
$${\Cal F}_+=\{F\in {\Cal F}| \mu_S(F) >0\}$$
and
$${\Cal F}_-=\{F\in {\Cal F}| \mu_S(F) <0\}.$$
Let $X,Y_1,\ldots, Y_m$ be variables corresponding to an  integral
basis of $\h^*$ such that
$$X|_{\s} \ \text{ corresponds \ to  \ an \ integral \ basis \ of } \ \s^* \
\text{ and } \ Y_j|_{\s}=0, \ j\geq 1.$$ We will identify
$$H_H^*(\text{pt})=S(\h^*)=\bC[X,Y_1, \ldots, Y_m].$$
The following  residue formula has been proved in  \cite{Je-Ki1}
and [Je-Ki2] (see also [Ka], [Gu-Ka]):

\proclaim{Theorem 1.2} For any $\eta\in H_H^*(N)$ we
have\footnote{Note that both sides of the equation are in
$\bC[Y_1, \ldots, Y_m]$.}\rm
$$\kappa_S(\eta)[N_{\text{red}}]= c\sum_{F\in {\Cal F}_+}\text{Res}_X^+
\frac{\eta|_F}{e_F}.\tag 1.1$$ \it Here $e_F\in H_H^*(F)$ denotes
the $H$-equivariant Euler class of the normal bundle $\nu(F)$, $c$
is a non-zero constant\footnote{The exact value of $c$
is not needed in our paper, but the interested reader can find it in
[Je-Ki2, \S 3].}, and \rm{$\text{Res}_X^+$} \it is defined
below in (1.3).
\endproclaim

According to [Je-Ki2, Definition 3.3 and Proposition 3.4], the value of
$\text{Res}_X^+$ on a rational function of the type
$$h(X,Y_1,\ldots,Y_m)=\frac{p(X,Y_1,\ldots, Y_m)}{\prod_k(m_kX+\sum_i\beta_{ki}Y_i)} \tag 1.2$$
where $p\in \bC[X,\{Y_i\}]$, $m_i\in \bZ\setminus\{0\}$ and $\beta_{ki}\in\bC$, can
be obtained
as follows:  regard $Y_1,\ldots,Y_m $ as constants (i.e. as  complex
numbers) and set
$$\text{Res}_X^+(h)=\sum_{b\in\bC}\text{Res}_{X=b}\frac{p}{q}\tag 1.3$$
where on the right hand side $\frac{p}{q}$ is interpreted as a
meromorphic function in the variable $X$ on $\bC$. It remains to
note that  only residues of expressions of the type (1.2)
 are involved in (1.1). For we may assume that each normal
bundle $\nu(F)$, $F\in{\Cal F}_+$, is a direct sum of
$H$-invariant line bundles and  in this way the  Euler class $e_F$
is the product of  $H$-equivariant first Chern classes of those line
bundles.

{\bf Remark 1.3.} For future reference, we mention that
$$\text{Res}_X^+\frac{1}{X+\sum_i\beta_iY_i}=1.$$

{\bf Remark 1.4.} Guillemin and Kalkman gave another definition of
the residue of  $h$ given by (1.2) (see  [Gu-Ka, section 3]).
Denote $\beta_k(Y)=\sum_i\beta_{ki}Y_i$ and write
$$\frac{1}{\prod_k(m_kX+\beta_k(Y))}=
\frac{1}{\prod_km_kX(1+\frac{\beta_k(Y)}{m_kX})}=\prod_k
(m_kX)^{-1}
\prod_k(1-\frac{\beta_k(Y)}{m_kX}+(\frac{\beta_k(Y)}{m_kX})^2
-\ldots ).$$ Multiply the right hand side by the polynomial $p$,
and add together all coefficients of $X^{-1}$: the result is the
Guillemin-Kalkman residue.  It is a simple exercise to show that
the latter coincides with the residue $\text{Res}_X^+$ of [Je-Ki2]
(one uses the fact that the sum of the residues of the 1-form $h(X)dX$ on the
Riemann sphere equals zero).

Since $\kappa_S$ is a ring homomorphism, from Theorem 1.2 we
deduce: \proclaim{Corollary 1.5}{\rm
$$\kappa_S(\eta)\kappa_S(\zeta)[N_{\text{red}}]=
\sum_{F \in {\Cal F}_+}\text{Res}_X^+
\frac{(\eta\zeta)|_F}{e_F}.$$}
\endproclaim

Since $0\in \s^*$ is a regular value of $\mu_S$,  $S$ acts on
$\mu_S^{-1}(0)$ with finite stabilizers. Hence $\mu_S^{-1}(0)/S$
has at worst orbifold singularities and in particular it satisfies
Poincar\'e duality. From Corollary 1.5 it follows that $\eta \in
H_H^*(N)$ is in $\ker \kappa_S$ iff
$$\sum_{F \in {\Cal F}_+}\text{ Res}_X^+
\frac{(\eta\zeta)|_F}{e_F}=0, \tag 1.4$$ for all $\zeta \in
H_H^*(N)$.

The main goal of our paper is to compare the latter description of
$\ker \kappa_S$ to the one given by Goldin in \cite{Go}. We will
assume  that the fixed point sets $N^S$ and $N^H$ are equal, which
holds for a generic choice of the subgroup $S$ of $H$ (see Remark
1.9 below).

 \proclaim{Theorem 1.6} {\rm (see [Go, Theorem 1.5])} We have
$$\ker \kappa_S =K_-\oplus K_+$$
where $K_-$ is the set of all elements of $H_H^*(N)$
vanishing on all components  $F\in {\Cal F}_-$, and similarly for $K_+$.
\endproclaim

The goal of our paper is to give a direct proof of the following
result:

\proclaim{Theorem 1.7} Assume that the fixed point set $N^H=N^S$
is finite. A class $\eta\in H_H^*(N)$ satisfies (1.4) for all
$\zeta \in H_H^*(N)$ if and only if $\eta \in K_+\oplus K_-$.
\endproclaim
  It is obvious that any $\eta \in K_+$ satisfies
(1.4). On the other hand, by the Atiyah-Bott-Berline-Vergne
localization formula (see \cite{At-Bo}, \cite{Be-Ve}) the sum
$$\sum_{F \in {\Cal F}}
\frac{(\eta\zeta)|_F}{e_F}$$ is a polynomial in $X, \{Y_i\}$, so that
its residue  is zero. It follows that any $\eta$  in $K_-$
satisfies (1.4) as well. The hard part will be to prove that if
$\eta$ satisfies (1.4) for any $\zeta\in H_H^*(N)$, then $\eta$ is
in $K_-\oplus K_+$.

{\bf Remark 1.8.}  Theorem 1.7 is a generalization of the main
result of \cite{Je}.

{\bf Remark 1.9.} Let $M$ be a symplectic manifold acted on by a
torus $G$ in a Hamiltonian fashion such that the fixed point set
$M^G$ is finite and let $\mu_G : M\to \g^*$ be the corresponding
moment map. Suppose that $0\in \g^*$ is a regular value of $\mu_G$
and consider the corresponding symplectic reduction
$$M\quott G=\mu_G^{-1}(0)/G$$  as well as the Kirwan surjection
$$\kappa: H_G^*(M)\to H^*(M\quott G).$$
Using reduction in stages, results like Theorem 1.6 or Theorem 1.7
can be used in order to obtain descriptions of $\ker\kappa$. The
following construction is needed: There exists a ``generic" (in
the sense of \cite{Go, section 1}) circle $S_1\subset G$ with Lie algebra
$\s_1$ such that

\item{$\bullet$} $M^{S_1}=M^G$

\item{$\bullet$} $0$ is a regular value of the moment map
$\mu_{S_1}:M\to \s_1^*$.

Consider the action of $G/S_1$ on the reduced space $M\quott
S_1=\mu_{S_1}^{-1}(0)/S_1$. As before, there exists a generic
circle $S_2\subset G/S_1$ such that

\item{$\bullet$} $(M\quott S_1)^{S_2}=(M\quott S_1)^{G/S_1}$

\item{$\bullet$} $0$ is a regular value of the moment map
$\mu_{S_2}:M\quott S_1 \to \s_2^*$.

The equivariant Kirwan map corresponding to the subtorus
$S_1\times S_2$ of $G$ is the following composition of maps:
$$H^*_{G}(M)@>{\kappa_1}>>H^*_{G/S_1}(M\quott S_1)@>{\kappa_2}>>
H^*_{G/(S_1\times S_2)}(M\quott (S_1\times S_2)),$$ where we have
used ``reduction in stages". In fact the procedure can be
continued giving rise to a sequence of tori
$$\{1\}=T_0\subset T_1\subset T_2 \subset \ldots \subset T_m=G$$
where $T_1=S_1, T_2=S_1\times S_2, \ldots $ with the following
properties:

\item{$\bullet$} $(M\quott T_{j-1})^{T_j/T_{j-1}}=(M\quott T_{j-1})^{G/T_{j-1}}$

\item{$\bullet$} $0$ is a regular value of the moment map
$\mu_{T_j/T_{j-1}}$ on $M\quott T_{j-1}$.

 By reduction in stages, the Kirwan map $\kappa$ decomposes as
$$\kappa=\kappa_{m}\circ \ldots \circ \kappa_{1} \tag 1.5$$
where
$$\kappa_j=\kappa_{T_j/T_{j-1}}:H_{G/T_{j-1}}^*(M//T_{j-1})
\to H^*_{G/T_j}(M//T_j)  $$ is the $T_j/T_{j-1}$-equivariant Kirwan
map.
 Goldin \cite{Go} used this decomposition in order to deduce from
 Theorem 1.6 the Tolman-Weitsman description of $\ker\kappa$ (see
 \cite{To-We}). Theorem 1.7 of our paper  --- with $N=M\quott T_{j-1}$,
 $H=G/T_{j-1}$ and $S=T_j/T_{j-1}$ --- shows that the same
 Tolman-Weitsman description
 of $\ker\kappa$ is equivalent to the obvious characterization in terms
 of residues arising from the components $\kappa_j$ of $\kappa$.

\head \S 2 Expressing  $\ker\kappa_S$ in terms of residues
\endhead

 We will give a  proof of Theorem 1.7.
Take $\xi\in \s$ a non-zero vector and
$f=\mu_S^{\xi}=\langle\mu_S,\xi \rangle :N\to \bR$ the  function
induced on $N$ by the moment map $\mu_S:N\to \s^*$. This is an
$H$-equivariant Morse function, whose critical set is  $N^S=N^H$.
For any $F\in{\Cal F}$, we have the $H$-equivariant splitting of
the tangent space
$$T_F N=\nu^-_f(F) \oplus \nu^+_f(F),$$
determined by the sign of the Hessian on the two summands. Let
$e(\nu^-_f(F)),e(\nu^+_f(F))\in H_H^*(F)$ be the corresponding
$H$-equivariant Euler classes. The following two results can be
proved by using Morse theory (see for example \cite{Go}):
\proclaim{Proposition 2.1}
 Suppose $\eta\in H_H^*(N)$ restricts to
zero on all $G\in{\Cal F}$ for which $f(G)<f(F)$. Then $\eta|_F$
is some \rm $H_H^*(\text{pt})$\it-multiple of $e(\nu^-_fF)$.
\endproclaim

\proclaim{Proposition 2.2} For any $F\in {\Cal F}$ there exists a
class $\alpha^-(F)\in H_H^*(N)$ with the following properties:

\item {1.} $\alpha^-(F)|_G=0$, for any $G\in {\Cal F}$ which cannot be
joined to $F$ along a sequence of integral lines of the negative
gradient field $-\nabla f$ (in particular, for any $G\in{\Cal F}$
with $f(G)<f(F)$)

\item {2.} $\alpha^-(F)|_F=e(\nu_f^-F)$.

In the same way there exists $\alpha^+(F)\in H_H^*(N)$ such that:

 \item {1.} $\alpha^+(F)|_G=0$, for any $G\in {\Cal F}$ which cannot be
joined to $F$ along a sequence of integral lines of the  gradient
field $\nabla f$ (in particular, for any $G\in{\Cal F}$ with
$f(G)>f(F)$)

\item {2.} $\alpha^+(F)|_F=e(\nu_f^+F)$.
\endproclaim

From the injectivity theorem of Kirwan
--- which says that $\alpha\in H_H^*(N)$ is uniquely determined by
its restrictions to  $N^H$ --- we deduce:

\proclaim{Corollary 2.3} The set $\{\alpha^-(F)|F\in {\Cal F}\}$
is a basis of  $H_H^*(N)$ as a $H_H^*(\text{pt})$-module.
\endproclaim

Consider the space $\H_H^*(N)$ consisting of all expressions of
the type
$$\sum_{F\in {\Cal F}} r_F\alpha^-(F),$$
where $r_F$ is in the ring $\bC(X,\{Y_i\})$ of rational expressions in
$X,\{Y_i\}$ (i.e. quotients $p/q$, with $p,q\in \bC[X,\{Y_i\}]$, $q\neq
0$). The space $\H_H^*(N)$ is obviously a $\bC(X,\{Y_i\})$-algebra.

\proclaim{Lemma 2.4} Take $\eta \in H_H^*(N)$ of degree $d$. Then
we can decompose
$$\eta=\eta_++\eta_-$$
where $\eta_+,\eta_-\in \H_H^*(N)$, such that

(i) $\eta_+|_{{\Cal F}_-}=0$,  $\eta_-|_{{\Cal F}_+}=0$;

(ii) $\eta_+$ and $\eta_-$ are linear combinations of
$\alpha^-(F)$, $F\in {\Cal F}$, where the coefficients $r_F$  are
rational functions whose denominators can be decomposed as
products of linear factors of the type
$$X+ \sum_i\beta_iY_i,\quad \beta_i\in \bC.$$
\endproclaim
\demo{Proof} By Corollary 2.3, we have
$$\eta=\sum_{F\in{\Cal F}}p_F\alpha^-(F),$$
where $p_F\in \bR[X,\{Y_i\}]$ are homogeneous polynomials. Take
$$\tilde{\eta}_-=\sum_{F\in{\Cal F}_-}p_F\alpha^-(F),\quad
\tilde{\eta}_+=\sum_{F\in{\Cal F}_+}p_F\alpha^-(F).$$ By
Proposition 2.2, $\tilde{\eta}_+$ restricts to 0 on all $G\in
{\Cal F}_-$.

Consider the ordering $F_1, F_2, \ldots $ of the elements of
${\Cal F}_+$ such that
$$0<f(F_1)<f(F_2)<\ldots .$$
There exists a rational function $r_1(X,\{Y_i\})\in \bC(X,\{Y_i\})$ such
that
$$\tilde{\eta}_-|_{F_1}=r_1(X,\{Y_i\})e(\nu_f^-F_1),$$
which means that the form
$$\tilde{\eta}_1:=\tilde{\eta}_--r_1(X,\{Y_i\})\alpha^-(F_1)$$
vanishes at $F_1$. There exists another rational function
$r_2(X,\{Y_i\})\in \bC(X,\{Y_i\})$, with the property that
$$\tilde{\eta}_1|_{F_2}=r_2(X,\{Y_i\})e(\nu_f^-F_2),$$
which implies that the form
$$ \tilde{\eta}_1-r_2(X,\{Y_i\})\alpha^-(F_2)$$
vanishes at both $F_1$ and $F_2$. We continue this process and we
get the decomposition claimed in the proposition as follows:
$$\eta_-=\tilde{\eta}_--r_1\alpha^-(F_1)-r_2\alpha^-(F_2)-\ldots, \quad
 \eta_+=\tilde{\eta}_++r_1\alpha^-(F_1)+r_2\alpha^-(F_2)+\ldots .
 $$
Property (ii) follows from the fact that for each $F\in{\Cal
F}_+$, the weights of the representation of $S$ on $\nu_f^-(F)$
are all non-zero.
 \qed
\enddemo

We are now ready to prove the main result of the paper.

\demo{Proof of Theorem 1.7} Take $\eta \in H_H^d(N)$ satisfying
$$
\text{Res}_X^+\sum_{F\in{\Cal F}_+}\frac{(\eta\zeta)|_F}{e_F}=0
\tag 2.1$$ for all $\zeta\in H_H^*(N)$. We consider the
decomposition
$$\eta=\eta_-+\eta_+$$
with $\eta_-, \eta_+\in \H_H^*(N)$ given by Lemma 2.4. We show that
$\eta_-$ and $\eta_+$ are actually in $H_H^*(N)$. More precisely,
if  $\eta_+$ is of the form
$$\eta_+=
\sum_{G\in{\Cal F}_+}\frac{p_G}{q_G}\alpha^-(G)
$$ with $p_G,q_G\in \bC[X,Y_i]$  we  show
that
$$q_G\ \text{divides} \ p_G$$
for any $G\in {\Cal F}_+$.

From (2.1) and the fact that $\eta_-|_{{\Cal F}_+}=0$, we deduce that
$$\text{Res}_{X}^+\sum_{F\in{\Cal F}_+}\frac{(\eta_+\zeta)|_F}{e_F}=0
$$
which is equivalent to
$$\text{Res}_{X}^+\sum_{F,G\in{\Cal F}_+}
\frac{p_G}{q_G}\cdot \frac{\alpha^-(G)|_F\zeta|_F}{e_F}=0 \tag
2.2$$ for all $\zeta \in H_H^*(N)$. Consider again the ordering
$F_1,F_2, \ldots $ of the elements of ${\Cal F}_+$ such that
$$f(F_1)<f(F_2)<\ldots .$$
We prove by induction on $k\geq 1$ that $q_{F_k}$ divides
$p_{F_k}$.

 In (2.2) we put
$\zeta=p\alpha^+(F_1)$, where $p\in \bC[X,Y_i]$ is an arbitrary
polynomial. Since
$$e(\nu_f^-F_1)e(\nu_f^+F_1)=e_{F_1},$$
we deduce that
$$\text{Res}_{X}^+p\frac{p_{F_1}}{q_{F_1}}=0$$
for any $p\in \bC[X,Y_1,\ldots,Y_m]$ (we are using the fact that
$\alpha^+(F_1)|_F = 0 $ if $\mu(F) > \mu(F_1) $ while
$\alpha^-(G)|_F = 0$ if $\mu(F) < \mu(G)$: hence the only nonzero
contribution comes from $F=G = F_1$). From Lemma 2.5 (see below)
we deduce that
$$q_{F_1}\ \text{divides} \ p_{F_1}.$$

Now we fix $k\geq 2$, assume that $q_{F_i}$ divides $p_{F_i}$ for
any $i<k$ and show that $q_{F_k}$ divides $p_{F_k}$. In (2.2) put
$\zeta =p\alpha^+(F_k)$, with $p\in \bC[X,Y_1,  \ldots,Y_m]$. We
obtain
$$\text{Res}_{X}^+\sum_{i\leq j \leq k}p\frac{p_{F_i}}{q_{F_i}}
\cdot \frac{\alpha^-(F_i)|_{F_j}\alpha^+(F_k)|_{F_j}}{e_{F_j}}=0.$$
The sum in the left hand side is over $i$ and $j$ ($k$ being fixed).
We divide into sums of the type
$$\Sigma^{i}:=\text{Res}_{X}^+\sum_{i\leq j \leq k}p\frac{p_{F_i}}{q_{F_i}}
\cdot \frac{\alpha^-(F_i)|_{F_j}\alpha^+(F_k)|_{F_j}}{e_{F_j}},\quad i\leq k.$$
 If $i< k$, then $\Sigma^{i}=0$, by the hypothesis
that $q_{F_i}$ divides $p_{F_i}$ and the
Atiyah-Bott-Berline-Vergne localization formula for
$p\frac{p_{F_i}}{q_{F_i}} \cdot \alpha^-(F_i)\alpha^+(F_k)$ (note
that $\alpha^-(F_i)|_F\alpha^+(F_k)|_F=0$ for any $F$ which is not
of the form $F_j$ with $i\leq j \leq k$).

Finally we show that
$\Sigma^{k}=0$, which is equivalent to
$$\text{Res}_{X}^+p\frac{p_{F_k}}{q_{F_k}}=0$$
for any $p\in \bC[X,Y_i]$. By Lemma 2.5,
$$q_{F_k}\ \text{divides} \ p_{F_k},$$
which concludes the proof.  \qed
\enddemo
We have used the following result: \proclaim{Lemma 2.5} Let $f,g$
be in $\bC[X,\{Y_i\}]$, where
$$g=\prod_k(X+\sum_i\beta_{ik}Y_i),$$ $\beta_{ik} \in \bC$. If \rm
$$\text{Res}_{X}^+(p\cdot\frac{f}{g})=0, \tag 2.3$$\it
for any $p\in \bC[X,\{Y_i\}]$, then $g$ divides $f$.
\endproclaim
\demo{Proof} Suppose that $g$ does not divide $f$. We can assume
that $g$ and $f$ are relatively prime. Then there exist
$p_1,p_2\in \bC[X,\{Y_i\}]$ such that
$$p_1f+p_2g=1.$$ From (2.3) it follows that
$$\text{Res}_{X}^+\frac{p}{g}=\text{Res}_{X}^+(pp_1\cdot\frac{f}{g}+p_2)=0$$
for any $p\in\bC[X,\{Y_i\}]$. Now fix $k_0$ and set
$$p=\prod_{k\neq k_0}(X+\sum_i\beta_{ik}Y_i).$$
We deduce that
$$\text{Res}_{X}^+\frac{1}{X+\sum_i\beta_{ik_0}Y_i}=0.$$
But the left hand side is actually 1 (see Remark 1.3), which is a
contradiction.\qed
\enddemo

\head \S 3 Residues and the kernel of the Kirwan map
\endhead

In this section we shall give a direct characterization of
the kernel of the Kirwan map in terms of residues of one variable.

Let $M$ be a compact symplectic manifold equipped with a Hamiltonian
$G$ action, where $G$ is a torus and let $\kappa:H_G^*(M)\to H^*(M\quott G)$
be the Kirwan map.
Consider the decomposition of $\kappa$ described in Remark 1.9
(see equation (1.5)). For each $j$ between $1$ and $m$ we consider the
commutative diagram
\newpage

$$
H^*_{G/T_{j-1}}(M\quott T_{j-1}) @>{\kappa_j}>> H^*_{G/T_j}(M \quott T_j)  $$
$$\downarrow \pi_j  \ \ \ \ \  \ \ \ \ \    \phantom{\longrightarrow } \downarrow \pi\tag 3.1$$
$$H^*_{T_j/T_{j-1}}(M\quott T_{j-1}) @>{\tilde{\kappa}_j}>> H^*(M \quott T_j)
$$
where $\tilde{\kappa}_j$ is the Kirwan map associated to the
action of the circle $T_j/T_{j-1}$ on $M\quott T_{j-1}$. For the
same action we consider the moment map $\mu_j$, the fixed point
set ${\Cal F}^j=(M\quott T_{j-1})^{T_j/T_{j-1}}$ and the partition
of the latter into ${\Cal F}^j_-$ and ${\Cal F}^j_+$ which consist
of fixed points where $\mu_j$ is negative, respectively
positive.

\proclaim{Definition 3.1} An element $\eta\in
H^*_{T_j/T_{j-1}}(M\quott T_{j-1})$ is in $\ker_{\text{res},j}$ if
 $$\text{Res}_{X_j = 0 } \sum_{F \in {\Cal F}^j_+ }
\frac{(\eta \zeta)|_{F} }{e_F(X_j)} = 0 $$ for any $\eta \in
H^*_{T_j/T_{j-1}}(M\quott T_{j-1})$, where $X_j$ is a variable
corresponding to a basis of the dual of the Lie algebra of
$T_j/T_{j-1}$ and $\text{Res}_{X_j=0}$ means the coefficient of
$X_j^{-1}$.\endproclaim

For any $\alpha\in H_G^*(M)$ we define
 $$\alpha_j = (\kappa_j \circ \dots \circ \kappa_1) (\alpha)$$
$0\leq j \leq m$.
The goal of this   section      is to prove

\proclaim{Theorem 3.2} Let $\alpha \in H^*_G(M)$. Then $\kappa
(\alpha) = 0 $ if and only if \rm $\pi_j( \alpha_{j-1} ) \in
\ker_{\text{res,j}}$ \it for some $j\geq 1$.
\endproclaim

\demo{Proof} If $\kappa(\alpha)=0$, by  equation (1.5) there
exists $j$ such that $\alpha_{j}=0$. Since
$\alpha_j=\kappa_j(\alpha_{j-1})$, we can use first the
commutativity of the diagram (3.1) to deduce that
$\tilde{\kappa}_j(\pi_j\alpha_{j-1})=0$ and then  the residue
formula of [Je-Ki1] and [Je-Ki2] to deduce that
 $\pi_j( \alpha_{j-1} ) \in \ker_{\text{res,j}}$.

The opposite direction is less obvious:
Take $\alpha\in H_G^*(M)$ of degree greater than zero such that
 $\pi_j \alpha_{j-1} \in \ker_{\text{res,j}}$.
If $j=m$ then $\pi_j$ and $\pi$ are the identity maps,
$\kappa_j=\tilde{\kappa}_j$ and we just have to apply the residue formula
of [Je-Ki1] and [Je-Ki2] to deduce that $\kappa_m(\alpha_{m-1})=
\kappa(\alpha)=0$. If $j<m$,   by the
same  residue formula, we have $\tilde{\kappa}_j(\pi_j \alpha_{j-1}) = 0 $,
which implies $\pi (\kappa_j \alpha_{j-1}) = 0 $.
Note that
$\ker(\pi) = H^*_{G/T_j} (\text{pt})$
so $ \pi(\kappa_j \alpha_{j-1}) = 0 $
 is equivalent to
$\kappa_j \alpha_{j-1} \in H^*_{G/T_j}({\text{pt}})$. The latter
implies that $(\kappa_m \circ \dots \circ \kappa_j) \alpha_{j-1} =
0 $ (since the map $\kappa_m \circ \dots \circ \kappa_{j+1} $
sends  all elements of $H^*_{G/T_j} ({\text{pt}})$ of degree
larger than zero  to $0$, because the image of this map is $H^*(M
\quott G)$ and the image of the equivariant cohomology of a point
under this map is the ordinary cohomology of a point). But this
means $\kappa(\alpha)=0$ and the proof is now complete. \qed
\enddemo

\bigskip

\eightpoint \it

\noindent Mathematics Department
\newline
University of Toronto
\newline
Toronto, Ontario M5S 3G3
\newline
Canada
\newline
{\tt jeffrey\@math.toronto.edu}

\bigskip

\it

\noindent Mathematics Department
\newline
University of Toronto
\newline
Toronto, Ontario M5S 3G3
\newline
Canada
\newline
 {\tt amare\@math.toronto.edu}

\enddocument